\documentclass[12pt]{amsart}
\usepackage{amssymb,amsmath,amsthm}

\begin{document}

\theoremstyle{plain}
\newtheorem{theorem}{Theorem}
\newtheorem{lemma}[theorem]{Lemma}
\newtheorem{proposition}[theorem]{Proposition}
\newtheorem{corollary}[theorem]{Corollary}
\newtheorem*{main}{Main Theorem}

\def\mod#1{{\ifmmode\text{\rm\ (mod~$#1$)}
\else\discretionary{}{}{\hbox{ }}\rm(mod~$#1$)\fi}}

\theoremstyle{definition}
\newtheorem*{definition}{Definition}

\theoremstyle{remark}
\newtheorem*{example}{Example}
\newtheorem*{remark}{Remark}
\newtheorem*{remarks}{Remarks}

\newcommand{\qq}{{\mathbb Q}}
\newcommand{\rr}{{\mathbb R}}
\newcommand{\nn}{{\mathbb N}}
\newcommand{\zz}{{\mathbb Z}}

\newcommand{\al}{\alpha}
\newcommand{\be}{\beta}

\newcommand{\rp}{\rr^+}
\newcommand{\ep}{\epsilon}
\newcommand{\lam}{\lambda}
\newcommand{\de}{\delta}

\thanks{This material is based in part upon work of the author,
  supported by the USAF under DARPA/AFOSR MURI Award
  F49620-02-1-0325. Any opinions, findings, and conclusions or
  recommendations expressed in this publication are those of the
  author and do not necessarily reflect the views of these agencies.}

\title{On the absence of uniform denominators in  Hilbert's 17th problem}

\author{Bruce Reznick}
\address{Department of Mathematics, University of 
Illinois at Urbana-Champaign, Urbana, IL 61801} 
\email{reznick@math.uiuc.edu}
\subjclass{Primary:  11E10, 11E25, 11E76, 12D15, 14P99}

\begin{abstract}
Hilbert showed that for most $(n,m)$ there exist psd forms
$p(x_1,\dots,x_n)$ of degree $m$ which cannot be written as
a sum of squares of forms. His 17th problem asked whether, in
this case, there exists a form $h$ so that $h^2p$ is a sum of
squares of forms; that is, $p$ is a sum of squares of rational
functions with denominator $h$. We show that, for every such
$(n,m)$ there does not exist a single form $h$ which serves
in this way as a denominator for {\it every} psd $p(x_1,\dots,x_n)$ of
degree $m$. 
\end{abstract}
\date{May 12, 2003}

\maketitle

\section{Introduction}

Let $H_d(\rr^n)$ denote the set of real homogeneous forms
of degree $d$ in $n$ variables (``$n$-ary $d$-ics'') .
By identifying $p \in H_d(\rr^n)$ with the $N = \binom
{n+d-1}{n-1}$-tuple of its coefficients, we
see that $H_d(\rr^n) \approx \rr^N$. Suppose $m$ is an even integer.
A form $p \in H_m(\rr^n)$ is called {\it positive semidefinite}
or {\it psd} if $p(x_1,\dots,x_n) \ge 0$ for all $(x_1,\dots,x_n) \in
\rr^n$. Following 
\cite{CL1}, we denote the set of psd forms in $H_m(\rr^n)$ by $P_{n,m}$.
Since $P_{n,m}$ is closed under addition and closed under
multiplication by positive scalars, it is a convex cone.
In fact, $P_{n,m}$ is a {\it closed} convex cone: if $p_n \to p$
coefficient-wise, and each $p_n$ is psd, then so is $p$. A psd form is called
{\it positive definite} or {\it pd} if $p(x_1,\dots,x_n) = 0$ implies
$x_j = 0$ for $1 \le j \le n$. The pd $n$-ary $m$-ics are the
interior of the cone $P_{n,m}$.

 A form $p \in H_m(\rr^n)$ is called  {\it a sum of squares} or
{\it sos}  if it can be written as a sum of 
squares of polynomials; that is, $p = \sum_k h_k^2$.
It is easy to show in this case that each $h_k \in H_{m/2}(\rr^n)$.
Again following \cite{CL1}, we denote the set
of sos forms in $H_m(\rr^n)$ by $\Sigma_{n,m}$. Clearly,
$\Sigma_{n,m}$ is a convex cone; less obviously, it is a
closed cone, a result due to R. M. Robinson \cite{Ro}.

In light of the inclusion  $\Sigma_{n,m} \subseteq P_{n,m}$, let
$\Delta_{n,m} = P_{n,m} \setminus \Sigma_{n,m}$. It was well-known by the
late 19th century that $P_{n,m} = \Sigma_{n,m}$ when $m=2$ or $n=2$.
In 1888, Hilbert proved \cite{Hi1} that $\Sigma_{3,4} = P_{3,4}$;
more specifically, 
 every $p \in P_{3,4}$ can be written as the sum of
three squares of quadratic forms. (An elementary proof,
with  ``five" squares is in \cite[pp.16-17]{CL2}; for modern expositions of
Hilbert's proof, see \cite{Sw} and \cite{Ru}.)  
Hilbert also proved in \cite{Hi1} that the preceding are the
{\it only} cases for 
which $\Delta_{n,m} = \emptyset$. That is, if $n\ge 3$ and $m \ge 6$
or $n \ge 4$ and $m \ge 4$, then there exist psd forms $n$-ary $m$-ics
that are not sos.  

In 1893, Hilbert \cite{Hi2} generalized his three-square result for
$P_{3,4}$ to ternary forms of higher degree.
Suppose $p \in P_{3,m}$ with $m \ge 6$. Then there exist 
$p_1 \in P_{3,m-4}$ and $h_{1k} \in H_{m-2}(\rr^3)$, $1 \le k \le 3$, so that
$$
 p_1p= h_{11}^2 + h_{12}^2 + h_{13}^2.
$$ 
(Hilbert's proof seems to be non-constructive, and lacks a modern
 exposition.  In the very recent paper  \cite{KP}, de Klerk and
Pasechnik discuss the implementation of an algorithm to find
  $p_1$ so that $p_1p$ is sos, though not necessarily as a
 sum of {\it three} squares. This paper uses Hilbert's result without
 giving an  independent proof.)

If $m = 6$ or 8, then $p_1$ is a sum of three squares of forms, and
hence (as Landau 
later noted \cite{L}), the four-square identity implies that $p_1^2p =
p_1(p_1p)$ is the sum of four squares of forms.
If $m \ge 10$, then the argument can be applied to $p_1$:
there exists $p_2 \in P_{3,m-8}$ with
$p_2p_1 = h_{21}^2 + h_{22}^2 + h_{23}^2$. Thus, if $m = 10$ or 12 (so
that $P_{3,m-8} = \Sigma_{3,m-8}$),
then $(p_1p_2)^2p = p_2(p_2p_1)(p_1p)$ is the sum of four squares of
forms,  An easy
induction shows that there exists $q \in H_t(\rr^3)$ with
$t = \lfloor \frac {(m-2)^2}8 \rfloor$ so that $q^2 p$ is the sum of four
squares of forms.

Hilbert's 17th Problem asked whether this generalizes to $n
> 3$ variables; that is, if $p \in P_{n,m}$, must there exist some
form $q$ so that $q^2p$ is sos?
Artin proved that there must be, in a way that gives no information
about $q$. Much more on the history of this subject can be found in the survey
paper \cite{Re2}.

This discussion leads to two closely related questions.
Suppose $p \in P_{n,m}$. Can we {\it find} a form $h$ such that $hp$ is
sos? Can we {\it find} a form $q$ so that $q^2p$ is sos?
If we've answered the second, we've answered the first. Conversely,
if $p\neq 0$ is psd and $hp$ is sos, then $h$ is psd. But it needn't be sos;
indeed, a trivial answer to the first question is to take $h=p$. 
 Stengle proved \cite{St} that if $p(x,y,z) =
x^3z^3 + (y^2z - x^3 - z^2x)^2$, then $p^{2s+1} \in
\Delta_{3,6(2s+1)}$ for every integer $s$. That is, $p^{2s-1}\cdot p$
is sos, but  $p^{2s}\cdot p$ is not.  Choi and Lam showed \cite{CL1}
that for $S \in \Delta_{3,6}$ (see (3) below), the product
$S(x,y,z)S(x,z,y)$ is actually sos.

The author gratefully acknowledges correspondence with Chip Delzell,
Pablo Parrilo, Vicki Powers, Marie-Fran\c{c}oise Roy and Claus
Scheiderer. Their suggestions have made this a better paper.
\section{What is known about the denominator}

The first concrete result about a denominator in Hilbert's 17th
Problem was found by P\'olya \cite{Po}. He showed that if $f \in
H_d(\rr^n)$ is positive on the unit simplex $\{(x_1,\dots,x_n)\ \vert\
x_j \ge
0, \sum x_j = 1\}$, then for sufficiently large $N$, $(\sum_j x_j)^Nf$
has positive coefficients. Replacing each $x_j$ by $x_j^2$, we see that if
$p \in H_{2d}(\rr^n)$ is  an even positive definite form, 
then $(\sum_j x_j^2)^Np$
is a sum of even monomials with positive coefficients, and so, as
it stands, is a sum of squares of monomials.  Taking even $N$, we see
that  $q = (\sum_j x_j^2)^{N/2}$ is a denominator for $p$.
 Habicht \cite{Ha} generalized P\'olya's proof to give an
alternate solution to Hilbert's 17th Problem for pd forms; however,
$h$ is not readily constructible and in general is no longer a power
of $\sum x_j^2$.  Except for one example, P\'olya did not
attempt to determine an explicit value of $N$. A good exposition of the
theorems of P\'olya and Habicht can be found in \cite{HLP}.

For positive definite  $p \in P_{n,m}$,  let 
$$\ep(p):= \frac {\inf \{
  p(u): u \in S^{n-1}\}}{\sup \{ p(u): u \in S^{n-1}\}}
$$
 measure how ``close" $p$ is to having a zero. The author \cite{Re1}
showed that if  
$$
N \ge \frac {nm(m-1)}{(4 \log 2)\ep(p)} -
\frac {n+m}2,
$$
then  $(\sum_j x_j^2)^Np$ is a sum
of $(m+2N)$-th powers of linear forms, and so is sos. A similar lower
bound has been shown to apply in P\'olya's case, one which goes to
infinity as $p$ approaches the boundary of $P_{n,m}$. (See papers by
de Loera and Santos \cite{LS}  and by Powers and the author \cite{PR1}.)

The restriction to positive definite forms is necessary.
There exist psd forms $p$ in $n \ge 4$ variables so that,
 if $h^2p$ is sos, then $h$ must have a specified zero.  
The existence of these unavoidable singularities, or  so-called ``bad
points",  insures that 
$(\sum x_j^2)^r p$ can never be a
sum of squares of forms for {\it any} $r$. Habicht's Theorem implies that
no positive definite
form can have a bad point. Bad points were first noted by  Straus
and have been extensively studied by Delzell; see, e.g. \cite{D1,D2}. 

\section{Recent results and a new theorem}

Scheiderer has shown in very recent work \cite{Sc} that for
$p \in P_{3,m}$, there 
exists $N = N(p)$ so that $(x^2 + y^2 + z^2)^N p(x,y,z)$ is sos;
indeed, $x^2+y^2+z^2$ can be replaced by any positive definite form.
This is a strong refutation to the existence of bad points  for
ternary forms.  

Also very recently, Lombardi and Roy \cite{LR} have constructed a
quantitative version of the Positivstellensatz. A special case is that
for fixed $(n,m)$, there exists $d = d(n,m)$ so that if $p \in
P_{n,m}$, there exists $q \in H_d(\rr^n)$ so that $q^2p$ is sos.

Suppose $(n,m)$ is such that $\Delta_{n,m} \neq \emptyset$. Theorem 1 below
states that there is no {\it  single} form $h$ so that, if $p \in P_{n,m}$,
then $hp$ is sos. Corollary 2 says that there is not even a {\it
  finite} set of forms $\mathcal H$ so that, if $p \in P_{n,m}$, 
 then there exists $h \in \mathcal H$ so that $hp$ in sos. In
 particular, there does not exist a 
finite set of denominators which apply to all of $P_{n,m}$.
This result implies that
$N(p)$ in Scheiderer's theorem is not bounded as $p$ ranges over
$P_{3,m}$. It also implies that the denominators in the Lombardi-Roy
theorem cannot be chosen from a finite, predetermined set. 

The proof of the Theorem is elementary and relies on a few simple
observations. If $p \neq 0$ is psd and $hp$ is sos, then $h$ is psd. 
As previously noted, $\Sigma_{n,m}$ is a closed cone for all
$(n,m)$. This cone is invariant under the action of taking invertible
linear changes of form. Thus, if $h'$ is derived from $h$ by such a
linear change, and if $hp$ is sos for every $p \in P_{n,m}$, then so
is $h'p$. Suppose
$\ell$ is a linear form, $p =  \sum_j g_k^2$ is sos, and $\ell \ |\
p$. Then $\ell^2 \ | \ p$ and $\ell \ | \ g_k$ for each $k$, and by
induction, $\ell^{2s} \ | \ p \implies \ell^s \ | \ g_k$. Thus, we
can ``peel off'' squares of linear factors from any sos form; this is
a common practice, dating back at least to \cite[p. 267]{Ro}.
 We use this observation in the contrapositive: if $p \in \Delta_{n,m}$,
then $\ell^{2s}p \in \Delta_{n,m+2s}$.

\begin{theorem} Suppose $\Delta_{n,m} \neq \emptyset$. Then there does
not exist a non-zero form $h$ so that if $p \in P_{n,m}$, then $hp$ is 
sos.
\end{theorem}
\begin{proof}
Suppose to the contrary that such a form $h$ exists.
 Since $h \neq 0$, there exists a point $a \in \rr^n$ so that $h(a)
 \neq 0$. By making an invertible linear change of variables, we can take
 $a = (1,0,\dots,0)$. Thus, we may assume without loss of generality
that $h(x_1,0,\dots,0) = \al x_1^d$, where $\al > 0$ and
$d$ is even. In the sequel, we 
distinguish $x_1$ from the other variables.

Choose $p \in P_{n,m} \setminus \Sigma_{n,m}$. Then
$$
h(x_1,x_2,\dots,x_n)p(x_1, r x_2,\dots, r x_n) 
$$
is sos for every $r \in \nn$. By making the change of variables $x_i \to
x_i/r$ for $i \ge 2$, we see that
$$
h(x_1,r^{-1}x_2, \dots, r^{-1}x_n)p(x_1,  x_2, \dots,  x_n) 
$$
is also sos. Since
$$
\lim_{r \to \infty} h(x_1,r^{-1}x_2, \dots, r^{-1}x_n) =
h(x_1,0,\dots,0) = \al x_1^d,
$$
and since $\Sigma_{n,m+d}$ is closed, it follows that
$$
\lim_{r \to \infty} h(x_1,r^{-1}x_2, \dots, r^{-1}x_n)p(x_1,
x_2, \dots,  x_n) = \al x_1^d p(x_1,\dots,x_n)
$$
 is sos. Thus $p$ is sos, a contradiction.
\end{proof}

The following elegant proof  is due to Claus Scheiderer and is
included with his permission; it supersedes the proof in an earlier
version of this manuscript.
\begin{corollary} Suppose $\Delta_{n,m} \neq \emptyset$. Then there does
not exist a finite set of non-zero forms $\mathcal H = \{h_1,...,h_N\}$ with
the property that, if $p \in P_{n,m}$, then $h_kp$ is 
sos for some $h_k \in \mathcal H$. 
\end{corollary}
\begin{proof}
Suppose $\mathcal H$ exists. For each $k$, there exists non-zero $p \in
\Delta_{n,m}$ so that $h_kp$ is sos. (Otherwise, we may delete $h_k$
harmlessly from  $\mathcal H$.)
Thus, each $h_k$ is psd, and there exists a form
$q_k$ so that $q_k^2h_k$ is sos. Define $h = \prod_k q_k^2 h_k$.
We now show that for every $p \in P_{n,m}$, $hp$ is sos: this 
contradicts the Theorem and proves the Corollary. By
hypothesis, there exists $h_j \in \mathcal H$ so that $h_j p$ is
sos. Thus,
$$
hp = \left( \prod_{k\neq j} q_k^2 h_k \right)\cdot q_j^2\cdot h_jp 
$$
is a product of sos factors, and so is sos.
\end{proof}

Finally, we know by Hilbert's theorem
that for $p \in P_{3,6}$, there exists quadratic $h$ so that $hp \in
\Sigma_{3,8}$. The three simplest forms in $\Delta_{3,6}$ are
\begin{equation}
M(x,y,z) = x^4y^2 + x^2y^4 + z^6 - 3x^2y^2z^2, \quad \text{due to
  Motzkin \cite{M}}; 
\end{equation}
Robinson's \cite{Ro} simplification of Hilbert's construction 
\begin{equation}
R(x,y,z) = x^6 + y^6 + z^6 - (x^4y^2 + x^2y^4 + x^4z^2 + x^2z^4 + y^4z^2 +
y^2z^4)  + 3x^2y^2z^2; 
\end{equation}
and
\begin{equation}
S(x,y,z) = x^4y^2 + y^4z^2 + z^4x^2 - 3x^2y^2z^2, \quad \text{due to
  Choi and Lam \cite{CL1,CL2}}.
\end{equation}

It is not too difficult to consider $qM, qR, qS$ for $q(x,y,z) = a^2x^2 +
b^2y^2 + c^2z^2$, and determine whether these are sos  using the
algorithm of \cite{CLR} directly  or its implementation in, e.g.,  \cite{Pa}.

Interestingly enough, these conditions are the same in each case: the forms
 are sos if and only if
$$
2(a^2b^2 + a^2c^2 + b^2c^2) \ge a^4 + b^4 + c^4.
$$
This expression factors rather neatly into:
$$
(a+b+c)(a+b-c)(b+c-a)(c+a-b) \ge 0,
$$
so if $a \ge b \ge c\ge 0$ without loss of generality, the only non-trivial
condition is that $b + c \ge a$; that is, there is a (possibly
degenerate) triangle with sides $a,b,c$. (Robinson \cite[p. 273]{Ro}
has a superficially similar condition, but note that his multiplier is
$a x^2 + b y^2 + c z^2$.)

By specializing this result and scaling variables as in the proof of
the theorem, we note that
$$
(x^2+y^2+z^2)M(x,\lambda y,\lambda z),\ (x^2+y^2+z^2)R(x,\lambda
y,\lambda z),\
(x^2+y^2+z^2)S(x,\lambda y,\lambda z)
$$
 are sos if and only if  $0 \le \lambda \le 2$.



\begin{thebibliography}{23}
\bibitem{CL1} Choi, M. D. and T. Y. Lam, \emph{An old question of
  Hilbert}, Queen's Papers in Pure and Appl. Math. (Proceedings of
Quadratic Forms Conference, Queen's University (G. Orzech ed.)), \textbf{46}
(1976), 385--405.

\bibitem{CL2}  Choi, M. D. and T. Y. Lam, \emph{Extremal positive
  semidefinite forms}, Math. Ann., \textbf{231} (1977), 1--18. 

\bibitem{CLR} Choi, M. D.,  T. Y. Lam and B. Reznick, \emph{Sums of
  squares of real polynomials}, Proc.  Sympos. Pure Math.,
  \textbf{58.2} (1995),  103--126. 

\bibitem{D1} Delzell, C. N., \emph{ Bad points for positive
  semidefinite polynomials}, Abstracts Amer. Math. Soc., \textbf{18}
  (1997),  \#926-12-174,  482.
 
\bibitem{D2} Delzell, C. N., \emph{Unavoidable singularities when
  writing polynomials as sums of squares of real rational functions},
  in preparation.

\bibitem{Ha} Habicht, W., \emph{ \"Uber die Zerlegung strikte
  definiter Formen in Quadrate}, Comment. Math. Helv., \textbf{12}
  (1940) 317--322. 

\bibitem{HLP} Hardy, G. H., J. E.. Littlewood and G. P\'olya,
  \emph{Inequalities}, Camb. U. Press, 2nd ed., 1967.

\bibitem{Hi1} Hilbert, D., \emph{ \"Uber die Darstellung definiter
  Formen als Summe von Formenquadraten}, Math. Ann. \textbf{32} (1888),
  342--350;   see Ges. Abh. 2, 154--161, Springer, Berlin, 1933,
  reprinted by Chelsea, New York, 1981.


\bibitem{Hi2} Hilbert, D., \emph{ \"Uber tern\"are definite Formen},
  Acta Math. \textbf{17} (1893) 169--197; see Ges. Abh. 2, 345--366, Springer,
  Berlin, 1933, reprinted by Chelsea, New York, 1981 .

\bibitem{KP} de Klerk, E. and D. V. Pasechnik, \emph{Products of
  positive forms, linear matrix inequalities, and Hilbert 17-th
  problem for ternary forms}, to appear in European J. of Oper.
  Res. 

\bibitem{L} Landau, E., \emph{\"Uber die Darstellung definiter
  Funktionen durch Quadrate}, Math. Ann., \textbf{62} (1906),
  pp. 272--285; also  in  Collected Works, vol. 2, pp.~237--250,
  Thales-Verlag, Essen, 1986.

\bibitem{LS} de Loera, J. A.  and F. Santos, \emph{ An effective
  version of P\'olya's theorem on positive definite forms},  J. Pure
  Appl. Algebra, \textbf{108} (1996), 231--240. (See correction, same
  journal, \textbf{155} (2001), 309--310.) 

\bibitem{LR} Lombardi, H. and M.-F. Roy, \emph{Elementary recursive
  degree bounds for Positivstellensatz}, in preparation.

\bibitem{M} Motzkin, T, S., \emph{The arithmetic-geometric
  inequality}, pp. 205--224  in \emph{Inequalities} (O. Shisha, ed.)
  Proc. of Sympos. at Wright-Patterson AFB, August 19--27, 1965 ,
  Academic Press, New York, 1967; also in  Theodore S. Motzkin:
  Selected Papers,  Birkh\"auser, Boston,  (D. Cantor, B. Gordon and
  B. Rothschild, eds.). 

\bibitem{Pa} Parrilo, P., \emph{Structured semidefinite programs and
  semialgebraic methods in robustness and optimization}, Ph.D. thesis,
 Calif. Inst. of Tech., 2000.

\bibitem{Po} P\'olya, G., \emph{ \"Uber positive Darstellung von
  Polynomen},  Vierteljschr. Naturforsch. Ges. Z\"urich, \textbf{73} (1928),
  141--145; see  Collected Papers, Vol. 2, pp. 309--313, MIT Press,
Cambridge, Mass., London, 1974.

\bibitem{PR1}  Powers, V. and B. Reznick, \emph{A new bound for
  P\'olya's theorem with applications to polynomials positive on
  polyhedra}, J. Pure Appl. Algebra \textbf{164} (2001), 221--229.


\bibitem{Re1} Reznick, B., \emph{Uniform denominators in Hilbert's
  Seventeenth Problem},   Math. Z., \textbf{220} (1995), 75--98.

\bibitem{Re2} Reznick, B., \emph{ Some concrete aspects of Hilbert's
  17th Problem},  Contemp. Math., \textbf{253} (2000), 251--272.

\bibitem{Ro} Robinson, R. M., \emph{Some definite polynomials which
  are not sums of squares of real polynomials}, Izdat. ``Nauka"
  Sibirsk. Otdel. Novosibirsk, (1973) pp. 264--282, (Selected questions
  of algebra and logic (a collection dedicated to the memory of
  A. I. Mal'cev), abstract in Notices AMS, \textbf{16} (1969), p. 554.

\bibitem{Ru} Rudin, W., \emph{Sums of squares of polynomials},
  Amer. Math. Monthly, \textbf{107} (2000), 813--821.

\bibitem{Sc} Scheiderer, C., \emph{Sums of squares on compact real
algebraic surfaces}, in preparation.

\bibitem{St} Stengle, G., \emph{Integral solution of Hilbert's
  seventeenth problem}, Math. Ann. \textbf{246} (1979/1980), 33--39.

\bibitem{Sw} Swan, R.G., \emph{Hilbert's theorem on positive ternary
  quartics}, Contemp. Math. \textbf{272} (2000), 287--292.
 


\end{thebibliography}
\end{document}